# Computing Zeta Functions Over Finite Fields

Daqing Wan

ABSTRACT. In this report, we discuss the problem of computing the zeta function of an algebraic variety defined over a finite field, with an emphasis on computing the reduction modulo $p^m$ of the zeta function of a hypersurface, where $p$ is the characteristic of the finite field.

## 1. Introduction

Let $p$ be a prime number. Let $\mathbb{F}_q$ be a finite field of $q$ elements of characteristic $p$. Let $X$ be an algebraic variety defined over $\mathbb{F}_q$, say an affine variety defined by the vanishing of $r$ polynomials in $n$ variables:

$$f_1(x_1, \cdots, x_n) = \cdots = f_r(x_1, \cdots, x_n) = 0,$$

where the polynomials $f_i$ have coefficients in $\mathbb{F}_q$. Let $N(X)$ denote the number of $\mathbb{F}_q$-rational points on $X$.

PROBLEM I. *Compute $N(X)$ efficiently.*

In addition to its intrinsic theoretical interest, this fundamental algorithmic problem has important applications in diverse areas such as coding theory, cryptography, primality testing, sphere packing, quasi-random number generator and fast multiplication. It is also very useful in bounding the number of torsion points of an abelian variety over a number field.

In principal, one can always check all the $q^n$ possibilities for $(x_1, \cdots, x_n)$. This is apparently very slow if either $q$ or $n$ is large. One would like to have faster non-trivial algorithms. In this paper, we discuss only deterministic algorithms in theoretical sense and we shall frequently use less precise but intuitive description.

For varieties with large automorphism groups (such as diagonal hypersurfaces), it is possible to find efficient elementary methods (such as Gauss sums) to compute $N(X)$. We shall, however, restrict to those methods which have potential to extend to general varieties.

In the special case that $X$ is an elliptic curve

$$E: \ y^2 = x^3 + ax + b, \ a, b \in \mathbb{F}_q,$$

a polynomial time algorithm was obtained by Schoof [Sc1]. More practical but probabilistic versions were obtained later by Atkin, Elkies, Couveignes and a few other authors, see [Sc2] for an updated exposition. Schoof's algorithm was generalized to abelian varieties and curves by Pila [Pi] with some improvements by Adleman-Huang [AH]. Curves and abelian varieties are the cases which have been studied most extensively in the literature. For more general varieties, no non-trivial algorithm is known unless the variety is defined over a small subfield of $\mathbb{F}_q$. It is shown in [GKS] that this problem is #P-hard (at least NP-hard) for sparse plane curves if one uses sparse input size. Thus, in order to get efficient algorithms to compute $N(X)$, one might need to put some restrictions on some of the parameters of $X/\mathbb{F}_q$, such as the dimension, degree, characteristic $p$, etc.

For each positive integer $k$, let $\mathbb{F}_{q^k}$ be the finite field of $q^k$ elements. This is the unique extension field of $\mathbb{F}_q$ of degree $k$. Let $N_k(X)$ denote the number of $\mathbb{F}_{q^k}$-rational points on $X$. We shall consider the following harder problem.

PROBLEM II. *Compute the sequence $\{N_k(X)\}_{k=1}^{\infty}$ efficiently.*

For fixed $n$, Problem I is easy if $q$ is small. One can simply check all $q^n$ possibilities for $x \in \mathbb{F}_q^n$. On the other hand, Problem II can already be hard even for small $q$ as $q^k$ will be large for large $k$. Nevertheless, one can still hope to find efficient algorithms in important cases that arise from various applications.

The sequence $N_k(X)$ has a nice structure. In fact, by Dwork's theorem [Dw], the generating zeta function

$$Z(X,T) = \exp(\sum_{k=1}^{\infty} \frac{N_k(X)}{k} T^k)$$

is a rational function. Thus, the sequence $N_k(X)$ satisfies a linear recurrence relation. Write

$$Z(X,T) = \frac{g(T)}{h(T)}$$

for some polynomials $g(T), h(T) \in \mathbb{Z}[T]$. Taking logarithmic derivative, one finds

$$\sum_{k=1}^{\infty} N_k(X) T^{k-1} = \frac{g'(T)}{g(T)} - \frac{h'(T)}{h(T)}.$$

Thus, Problem II is equivalent to

PROBLEM III. *Compute the zeta function $Z(X,T)$ efficiently.*

If an efficient algorithm to compute $N(X)$ is known, then one can continue to compute $N_k(X)$ for all positive integers $k$ up to something like the "degree" of the zeta function $Z(X,T)$. This easily leads to what I called the naive algorithm to compute $Z(X,T)$. It works well if the degree of $Z(X,T)$ is small and if one knows how to compute $N(X)$ efficiently. In particular, one obtains a good algorithm if the degree of $Z(X,T)$ is small and if $X$ is defined over a small subfield of $\mathbb{F}_q$ ($q$ can be large). But the degree of the zeta function can be very large. For a smooth projective hypersurface of dimension $n$ and large degree $d$, the degree of the zeta function is about the size $d^{n+1}$.

For an elliptic curve $E$, the zeta function is determined by $N_1(E)$:

$$Z(E,T) = \frac{1 + c(E)T + qT^2}{(1-T)(1-qT)},$$

where $c(E) = N_1(E) - (q+1)$. Thus, the zeta function $Z(E,T)$ can be computed efficiently by Schoof's work. More generally, when $X$ is an abelian variety of a fixed dimension embedded in a fixed projective space, the zeta function $Z(X,T)$ can also be computed efficiently by the work of Pila and Adleman-Huang.

For a non-singular plane curve of large degree $d$ over a small finite field $\mathbb{F}_q$ ($q = 2$ or 3 for instance), computing $N_1(X)$ is easy since $q$ is small. Even in this special case, computing $Z(X,T)$ already seems to be difficult as Katz-Sarnak pointed out [Po]. This case was motivated by their investigation of the distribution of the zeros and poles of zeta functions. The results of Adleman-Huang give a doubly exponential algorithm for a general plane curve, although they have a significantly faster algorithm for hyperelliptic curves (still exponential in terms of the degree). This is far from being practical if the degree is large.

Computing the zeta function of a hyperellitptic curve is very useful in the hyperelliptic curve cryptosystem proposed by Koblitz [Ko1-2]. Additional applications may be found in [AH2], [AG] and [Ts]. Computing zeta functions of algebraic varieties should also be the first key step in computing zeta functions of more general Hilbert sets and definable sets, see [Wa1] and [FHJ]. A very special case is already considered in [GKS].

Our purpose here is to give a brief general discussion of the modular approach for computing zeta functions, indicating some theoretical tools that are available and some realistic results that one might hope. As supporting evidence, we shall also describe some elementary theorems on computing the reduction of $Z(X,T)$ modulo $p^m$.

## 2. Modular approach for computing zeta functions

Let $X$ be an algebraic variety defined over a finite field $\mathbb{F}_q$. The degree of $Z(X,T)$ can be estimated explicitly using $p$-adic methods as done by Bombieri [Bo]. The size of the coefficients in $Z(X,T)$ can also be bounded easily using trivial estimate. In nicer cases, one can use Deligne's deep result [De] to get better bounds. The general modulo idea to compute $Z(X,T)$ consists of two steps. The first step is to compute $Z(X,T)$ modulo various small primes (or prime powers). The second step is to use the Chinese remainder theorem to recover $Z(X,T)$ from its various reductions. The second step is pretty easy and standard. The difficulty lies in the first step, namely, computing the reduction of $Z(X,T)$ modulo various small primes or prime powers.

Zeta functions of higher dimensional varieties over finite fields have been studied extensively from theoretic point of view, motivated by Weil's conjectures. Several theories are available but all of them are highly non-trivial and fairly difficult. They can be broadly classified as $\ell$-adic methods and $p$-adic methods, where $\ell$ is any prime number different from $p$. In this section, we give a brief speculative discussion and perspective about the potential of these methods in computing zeta functions. No results and/or theorems are given in this section.

**$\ell$-adic methods.** One can try to use Grothendieck's $\ell$-adic trace formula [Gr] modulo $\ell^m$:

$$Z(X,T) \equiv \prod_{i=0}^{2\dim(X)} \det(I - FT | H_c^i(X \otimes_{\mathbb{F}_q} \bar{\mathbb{F}}_q, \mathbb{Z}/\ell^m\mathbb{Z})^{(-1)^{i-1}} \pmod{\ell^m},$$

where $H_c^i$ denotes the $\ell$-adic étale cohomology with compact support and $F$ is the geometric Frobenius map acting on $H_c^i$. Unfortunately, in the general case, one does not know how to explicitly construct the $\ell$-adic cohomology and its Frobenius map. In fact, for singular or open varieties, one does not even know any explicit bound for the dimension of the $\ell$-adic cohomology $H_c^i$. Thus, a major future research problem is to make the $\ell$-adic theory constructive and explicit. If successful, it would have dramatic consequence on the problem of counting the number of rational points of a variety over a finite field.

In the special case of elliptic curves and abelian varieties, the $\ell$-adic cohomology can be constructed quite explicitly using torsion points and division polynomials. This immediately leads to the algorithms of Schoof for elliptic curves and Pila for abelian varieties. For a smooth curve $X$ of genus $g$, one can use its Jacobian variety (an abelian variety) and Pila's algorithm to compute the zeta function $Z(X,T)$ of the curve $X$, as done by Adleman-Huang. This algorithm is in general at least doubly exponential in $g$. In the case of hyperelliptic curve, the running time has been improved to be exponential in $g$. Thus it is still a very slow algorithm if $g$ is large. If $g$ is very small, this gives a good polynomial time algorithm. To conclude, even for smooth plane curves of large degrees over a small finite field, there is still no efficient algorithm to compute the zeta function. If $X$ is singular, there is also a problem of resolving the singularities if one wants to use this approach.

**$p$-adic methods**. There are various $p$-adic formulas for $Z(X,T)$. All of them can be made to be explicit in some sense. The earliest one is due to Dwork, with generalizations by Reich-Monsky [Mo]. There are also $p$-adic cohomological formulas in terms of formal cohomology, crystalline cohomology and rigid cohomology [Ber]. We believe that some of these formulas and their variants can be very useful in computing the zeta function in the general case if the characteristic $p$ is not too large. In particular, we expect it to be plausible to construct a $p$-adic algorithm to compute $Z(X,T)$, which runs in polynomial time (in the input and output size) if $X$ is sparse and the characteristic $p$ is small ($q$ can be large). This is exactly the case that is most useful in coding theory and cryptography. A weak corollary would be a polynomial time algorithm to compute $N(X)$ for $X$ of low degrees over a large finite field of small characteristic. Note that an elliptic curve is sparse and of low degree. If $p$ is very large, we do not even know how to efficiently compute the reduction $Z(X,T)$ modulo $p$.

It is worthwhile to mention that the $p$-adic étale cohomological approach is not suitable for computing the whole zeta function. This is because the $p$-adic étale formula gives only the $p$-adic unit root part of the zeta function, not the whole thing, as conjectured by Katz [Ka1] and proved in the zeta case by Etesse-Le Stum [ES].

Our general feeling is as follow. If $p$ is small ($q$ can be large), $p$-adic methods should be much more efficient. If $p$ is large and the degree of $X$ is small, then $\ell$-adic methods should be better assuming that the $\ell$-adic cohomology can be constructed explicitly. If both $p$ and the degree are large, it is not clear if the $p$-adic methods or the $\ell$-adic methods would work very well. A combination of both $p$-adic methods and $\ell$-adic methods could be useful in certain situations.

For $p$-adic methods, the major advantage is that everything can be made to be explicit. This at least provides explicit algorithms which work in the general case and which are faster in many cases than the trivial counting method. Another possible advantage is that $p$-adic methods easily extend to L-functions of exponential sums with the same complexity bound. The drawback is that one does not have any flexibility to choose any other (small) prime to work with ( unlike the $\ell$-adic case). Thus, if $p$ is large, one is stuck and even the modulo $p$ information is already something quite substantial.

For $\ell$-adic methods, the major advantage is that one can always choose many small primes $\ell \neq p$ to work with. One expects that it would be easier and more efficient to compute $Z(X, T)$ modulo a small prime. However, the $\ell$-adic cohomology is not explicit yet and thus one does not even know how to compute $Z(X, T)$ using $\ell$-adic methods in general. Even if the $\ell$-adic methods are made explicit, the resulting algorithm could be very slow if the degree of the zeta function is large. For L-functions of exponential sums over a finite field of large characteristic $p$, it is not clear if the $\ell$-adic methods would work well. This corresponds to the case of zeta functions of large degrees.

In the rest of this paper, we discuss the weaker problem of computing the reduction of zeta functions modulo $p$ or a power of $p$. In section 3, we consider the geometrically trivial case of zero dimensional hypersurfaces. In section 4, we consider the general problem of computing the zeta function modulo $p$ for higher dimensional varieties focusing on hypersurfaces. Finally, in section 5, we compute $Z(X, T)$ modulo a higher but small power of $p$.

## 3. Zero dimensional hypersurfaces

Let $X$ be the zero dimensional hypersurface defined by a polynomial $f(x) \in \mathbb{F}_q[x]$ in one variable of degree $d$. Geometrically, $X$ is trivial. It is a disjoint union of several closed points. Computing the zeta function in this simplest case is already non-trivial although it can be done efficiently as we shall see below.

Let
$$f(x) = \prod_{i=1}^{k} f_i(x)^{a_i}, \ a_i \geq 1$$
be the standard factorization of $f(x)$ over $\mathbb{F}_q$. Let $d_i$ $(1 \leq i \leq k)$ be the degree of the irreducible factor $f_i(x)$. An irreducible factor of $f(x)$ of degree $d_i$ corresponds to a closed point of $X$ of degree $d_i$. The Euler product form of $Z(X, T)$ shows that

$$Z(X, T) = \prod_{i=1}^{k} \frac{1}{1 - T^{d_i}}. \tag{3.1}$$

QUESTION 3.1. *Compute the zeta function $Z(X, T)$ directly from $f(x)$ without factoring $f(x)$ over $\mathbb{F}_q$.*

By formula (3.1), an algorithm for factoring $f(x)$ easily gives an algorithm for computing $Z(X, T)$. It is interesting to note that the converse seems also true quite often. In this section, we shall discuss three elementary formulas for $Z(X, T)$ modulo $p$, related to each other by duality and conjugacy. These lead to three simple polynomial time algorithms to compute $Z(X, T)$. Each formula leads to an efficient algorithm for factoring $f(x)$ if $p$ is small. One formula leads to the classical algorithm of Berlekamp. Another formula leads to the more recent algorithm of Niederreiter. The third one leads to a new algorithm, which has the advantage to extend to obtain higher dimensional zeta functions modulo $p$.

### 3.1. Zeta function and Berlekamp's algorithm

Let $R = \mathbb{F}_q[x]/(f)$ be the coordinate ring of the variety $X/\mathbb{F}_q$. This is an $\mathbb{F}_q$-vector space of dimension $d$. Let $F$ be the $q$-th power Frobenius map on $\mathbb{F}_q[x]$:

$$F : h(x) \longrightarrow h(x)^q.$$

It fixes the ideal $(f)$ and thus induces a map on $R$. Using the Chinese remainder theorem and normal basis, one can easily prove the following result due to Petr (1937).

THEOREM 3.1.1. *We have the congruence formula*

$$Z(X,T)^{-1} \equiv \det(I - FT|R) \pmod{p}.$$

The proof immediately gives the following result of Butler (1954):

THEOREM 3.1.2. *Let $R_1(F)$ denote the subspace of $R$ which is fixed by $F$. Namely, $R_1(F)$ is the eigenspace of $F$ at the eigenvalue $1$. Then*

$$\dim_{\mathbb{F}_q} R_1(F) = k,$$

*where $k$ is the number of distinct irreducible factors of $f(x)$ over $\mathbb{F}_q$.*

Applying this result to $f(x)$ over the extension field $\mathbb{F}_{q^j}$, we obtain the following result of Schwarz (1956).

THEOREM 3.1.3. *Let $k(j)$ be the number of distinct irreducible factors of $f(x)$ over $\mathbb{F}_{q^j}$. Let $s_i$ be the number of distinct irreducible factors of degree $i$ of $f(x)$ over the ground field $\mathbb{F}_q$. Then, we have*

$$\dim_{\mathbb{F}_q} R_1(F^j) = k(j)$$
$$= \sum_{i=1}^{d} (i,j) s_i,$$

*where $(i,j)$ denotes the g.c.d. of the positive integers $i$ and $j$.*

Let $A$ be the $d \times d$ matrix whose $ij$ entry is $(i,j)$, where $1 \leq i,j \leq d$. As indicated by Schwarz, the matrix $A$ is invertible. In fact, an explicit formula for the inverse matrix of $A$ can also be found, see [GA].

The map $F$ acting on $R$ can be computed efficiently by repeated squaring. The above theorems show that all the $k(i)$ and hence all the $s_i$ ($1 \leq i \leq d$) can be computed in polynomial time. Thus, the zeta function

$$Z(X,T) = \prod_{i=1}^{d} \frac{1}{(1-T^i)^{s_i}}$$

can be computed in polynomial time. It turns out that the subspace $R_1(F)$ is also useful for the harder problem of factoring $f(x)$ into irreducible factors over $\mathbb{F}_q$.

THEOREM 3.1.4. *(Berlekamp [Be]).* The eigenspace $R_1(F)$ can be used to factor $f(x)$ over $\mathbb{F}_q$.

The idea is as follows. By the Chinese remainder theorem, the eigenspace $R_1(F)$ has a basis of the form

$$\{y_1(x)\frac{f(x)}{f_1(x)^{a_1}}, \cdots, y_k(x)\frac{f(x)}{f_k(x)^{a_k}}\}$$

where $y_i(x) \in R$ satisfies

$$y_i(x)\frac{f(x)}{f_i(x)^{a_i}} \equiv 1 \ (\mod f_i(x)^{a_i}).$$

This basis shows that if $h_1(x)$ and $h_2(x)$ are two elements of $R_1(F)$ which are linearly independent over $\mathbb{F}_q$, then

$$f(x) = \prod_{c \in \mathbb{F}_q} (f(x), h_1(x) - ch_2(x)) \qquad (3.1.1)$$

gives a non-trivial factorization of $f(x)$. This algorithm clearly runs in polynomial time if $q$ is small. It can be extended to the case when $p$ is small ($q$ can be large). At present, no deterministic polynomial time algorithm is known to factor $f(x)$ in general if $p$ is large even under the assumption of the Generalized Riemann Hypothesis.

A subspace $R_1$ of $R$ is called **admissible** if any two linearly independent elements $h_1$ and $h_2$ of $R_1$ satisfy (3.1.1). It is clear that if an admissible subspace of $R$ can be constructed efficiently, then one can proceeds in the same way to construct a new factorization algorithm. In the next two subsections, we describe two more such admissible subspaces and hence two more such algorithms for factoring polynomials. Admissible subspaces are studied in some details in Lee-Vanstone [LV]. In particular, the congruence formula in section 3.3 gives a new admissible subspace and hence answers a question in [LV].

### 3.2. Zeta function and Niederreiter's algorithm

Let $\psi_q$ be the $\mathbb{F}_q$-linear operator on $\mathbb{F}_q[x]$ defined by

$$\psi_q(x^u) = \begin{cases} x^{u/q}, & \text{if } q|u, \\ 0, & \text{otherwise.} \end{cases}$$

This is a one sided inverse of the Frobenius map. It arises in various contexts such as Dwork's operator in $p$-adic theory of zeta functions, Cartier's operator in algebraic geometry and Hecke's U-operator in modular forms.

Let $H^{(q-1)}$ denote the $(q-1)$-th Hasse derivative on $\mathbb{F}_q[x]$ defined by

$$H^{(q-1)}(x^u) = \binom{u}{q-1} x^{u-(q-1)}.$$

Let $D$ be the $\mathbb{F}_q$-linear differential operator on $\mathbb{F}_q[x]$:

$$D = \psi_q \circ H^{(q-1)} \circ f^{(q-1)},$$

where $f^{q-1}$ is the multiplication operator. Since

$$\psi \circ f^q = f \circ \psi_q, \ f^q \circ H^{(q-1)} = H^{(q-1)} \circ f^q,$$

the ideal $(f)$ of $\mathbb{F}_q[x]$ is stable under the action of $D$. Thus, the operator $D$ induces a well defined map on $R = \mathbb{F}_q[x]/(f)$. The subspace $R_1(D)$ of $R$ fixed by $D$ has the admissible basis
$$\{f_1'(x)\frac{f(x)}{f_1(x)}, \cdots, f_k'(x)\frac{f(x)}{f_k(x)}\},$$
where $f_i'(x)$ denotes the derivative of $f(x)$. It follows that we have

THEOREM 3.2.1. *(Niederreiter [Ni]). The eigenspace $R_1(D)$ can be used to factor $f(x)$ over $\mathbb{F}_q$ in a way similar to Berlekamp's algorithm.*

Using the Chinese remainder theorem, normal basis and basic properties of the operators $\psi_q$ and $H^{(q-1)}$, it is elementary to prove directly

THEOREM 3.2.2. *We have the congruence formula*
$$Z(X,T)^{-1} \equiv \det(I - DT|R) \pmod{p}.$$

As before, the eigenspaces $R_1(D^j)$ can be used to determine the integers $s_i$ and hence gives a polynomial time algorithm to compute the zeta function $Z(X,T)$.

### 3.3. A new congruence formula

Let $G$ be the composition of $\psi_q$ and the multiplication operator $f^{q-1}$ acting on $\mathbb{F}_q[x]$:
$$G = \psi_q \circ f^{(q-1)}.$$
Again, the ideal $(f)$ of $\mathbb{F}_q[x]$ is stable under the action of $G$. Thus, the operator $G$ induces a well defined map on $R = \mathbb{F}_q[x]/(f)$. Using the Chinese remainder theorem, normal basis and basic property of the operators $\psi_q$, we deduce

THEOREM 3.3.1. *Assume that $f(0) \neq 0$. Then,*
$$Z(X,T)^{-1} \equiv \det(I - GT|R) \pmod{p}.$$

Note that this result is false in general if $f(0) = 0$. Similarly, the eigenspaces $R_1(G^j)$ can be used to determine the integers $s_i$ and hence yields a polynomial time algorithm to compute the zeta function $Z(X,T)$. For $f(0) \neq 0$, the eigenspace $R_1(G)$ has the admissible basis
$$\{xf_1'(x)\frac{f(x)}{f_1(x)}, \cdots, xf_k'(x)\frac{f(x)}{f_k(x)}\}.$$
Thus, we obtain

THEOREM 3.3.2. *The eigenspace $R_1(G)$ can be used to factor $f(x)$ over $\mathbb{F}_q$ in a way similar to the algorithms of Berlekamp and Niederreiter.*

### 3.4. Duality and conjugacy

The above three operators are closely related to each other. Using the definitions of $D$ and $G$, one can easily check that they are conjugate to each other via the element $x$ if $x$ is invertible. Namely, we have

THEOREM 3.4.1. *If $f(0) \neq 0$, then $D = x^{-1}Gx$.*

The relation to the Frobenius map is given by duality. The proof is elementary.

THEOREM 3.4.2. *There are two non-singular pairings on $R \times R$ such that $F$ and $D$ (resp. $F$ and $G$) are dual with respect to the first (resp. the second) pairing.*

The advantage of the operator $E = \psi_q \circ f^{q-1}$ is that it extends to higher dimensional varieties. We shall discuss this in next section.

## 4. Higher dimensional varieties

For simplicity, we restrict to the case of a hypersurface. The method works for an arbitrary variety. Thus, let $X$ be the affine hypersurface defined by a polynomial $f(x_1, \cdots, x_n)$ over $\mathbb{F}_q$ of total degree $d$ in $n$ variables. Let $R(d)$ be the $\mathbb{F}_q$-vector space
$$R(d) = (x_1 \cdots x_n \mathbb{F}_q[x_1, \cdots, x_n])_{\leq d}.$$
Namely, $R(d)$ is the $\mathbb{F}_q$-vector space generated by those monomials of degree at most $d$ and divisible by the product $x_1 \cdots x_n$. One computes that
$$\dim_{\mathbb{F}_q} R(d) = \binom{d}{n}.$$
The operator $\psi_q$ extends easily to several variable case. We shall see that the operator $\psi_q \circ f^{q-1}$ can be used to compute the zeta function $Z(X, T)$ modulo $p$.

PROPOSITION 4.1. *The operator $\psi_q \circ f^{q-1}$ is stable on the subspace $R(d)$.*

**Proof**. Let $h \in R(d)$. Then, $h$ is a polynomial of degree at most $d$. One checks that the degree of $f^{q-1}h$ is at most $d(q-1) + d = dq$. Thus, the degree of $\psi_q(f^{q-1}h)$ is at most $d$. Furthermore, if $h$ is divisible by $x_1 \cdots x_n$, then $\psi_q(f^{q-1}h)$ is also divisible by $x_1 \cdots x_n$. This proves that $\psi_q(f^{q-1}h) \in R(d)$.

The following simple congruence formula for zeta functions can be easily proved using the reduction of the Dwork trace formula. It is implicit in section 7 of [Wa2].

THEOREM 4.2. *We have the congruence formula*
$$Z(X, T)^{(-1)^n} \equiv \det(I - (\psi_q \circ f^{q-1})T | R(d)) \ (\mod p).$$

In the projective case, a similar congruence formula but with a much harder proof is given by Katz [Ka2]. Very general congruence formulas but acting on infinite dimensional space can be found in [TW1] and [TW2]. Using Theorem 4.2, we get

COROLLARY 4.3. *The zeta function $Z(X, T)$ modulo $p$ can be computed in time that is a polynomial in $\binom{d}{n} \log q$ (the input and output size) if $p$ is small and $n$ is fixed.*

This follows directly from Theorem 4.2 if $q$ is small. If $q$ is large but $p$ is small, one needs to use a little semi-linear algebra and Galois theory. If $p$ is large, we do not get a polynomial time algorithm. In computing the high power $f^{q-1}$, we have to expand it without reduction and the size simply gets larger and larger if $p$ is large.

## 5. Higher power congruences

The congruence result in section 4 can be improved to get finer information by using congruences modulo a higher power of $p$. For this purpose, we need to use $p$-adic liftings. Let $K$ be an unramified finite extension of $\mathbb{Q}_p$ with residue field $\mathbb{F}_q$. Let $\mathcal{O}_K$ be the ring of integers in $K$. For a positive integer $m$, Let $\mathcal{O}_m$ be the residue class ring $\mathcal{O}_K/(p^m)$, which is a finite ring of $q^m$ elements.

We still restrict to the case of a hypersurface. For simplicity of description, we shall only consider those solutions with $x_1 \cdots x_n \neq 0$. Thus, let $X$ be the affine hypersurface in the torus $\mathbb{G}_m^n/\mathbb{F}_q$ defined by a polynomial $f(x_1, \cdots, x_n)$ over $\mathbb{F}_q$ of total degree $d$ in $n$ variables. We assume that the polynomial $f$ has already been lifted to $\mathcal{O}_m$ for some $m > 0$. We shall consider the question of computing $Z(X,T)$ modulo $p^m$.

Let $R_{m,d}$ be the free $\mathcal{O}_m$-module:

$$R_{m,d} = \mathcal{O}_m[x_1, \cdots, x_n]_{\leq dp^{m-1}}.$$

Namely, $R_{m,d}$ is the free $\mathcal{O}_m$-module generated by those monomials of degree at most $dp^{m-1}$. One computes that

$$\dim_{\mathcal{O}_m} R_{m,d} = \binom{dp^{m-1}+n}{n}.$$

The operator $\psi_q$ is defined as before. It is easy to check that we have

PROPOSITION 5.1. *The operator $\psi_q \circ f^{(q-1)p^{m-1}}$ is stable on the finite dimensional free $\mathcal{O}_m$-module $R_{m,d}$.*

The operator $\psi_q \circ f^{(q-1)p^{m-1}}$ modulo $p^m$ is clearly independent of the choice of the lifting of $f$ to $\mathcal{O}_m$. The following simple congruence formula for zeta functions can be proved using the reduction modulo $p^m$ of the Dwork trace formula.

THEOREM 5.2. *We have the congruence formula*

$$\left(\frac{Z(X,T)}{Z(\mathbb{G}_m^n,T)}\right)^{(-1)^n} \equiv \prod_{i=0}^n \det(I - q^i(\psi_q \circ f^{(q-1)p^{m-1}})T | R_{m,d})^{(-1)^i \binom{n}{i}} \pmod{p^m}.$$

Since the zeta function of the $n$-torus $\mathbb{G}_m^n$ is very simple, we are reduced to computing $\det(I - (\psi_q \circ f^{(q-1)p^{m-1}})T | R_{m,d})$ modulo $p^m$. Using the above formula, we can show

COROLLARY 5.3. *$Z(X,T)$ modulo $p^m$ can be computed in polynomial time if $p^m$ is small and $n$ is fixed.*

This follows from Theorem 5.2 if $q$ is also small. If $q$ is large but $p^m$ is small, it involves a little extra work. If $p^m$ is large, we do not get a polynomial time algorithm as the dimension of $R_{m,d}$ is already exponential.

The approach we take in this section is the most elementary one. It has the advantage to be a polynomial time algorithm for any hypersurface as long as $p^m$ is small and as long as we only want modulo $p^m$ information for $Z(X,T)$. We shall consider another deeper $p$-adic approach in a future article, which seem to work well in some interesting cases even when we want the full information for $Z(X,T)$.

## References

[AH]   L.M. Adleman and M.D. Huang, Counting rational points on curves and abelian varieties over finite fields, Lecture Notes in Computer Science, 1122 (1996), 1-16.

[AG]   M. Anshel and D. Goldfeld, Zeta functions, one-way functions and pseudorandom number generators, Duke, Math., J., 88(1997), no.2, 371-390.

[Be]   E.R. Berlekamp, Factoring polynomials over finite fields, Bell System Technical J., 46(1967), 1853-1859.

[Ber]  P. Berthelot, Géométrie rigide et cohomologie des variétés algébriques de caractéristique $p$, Mém. Soc. Math. France (N.S.) No. 23(1986), 3, 7–32.

[Bo]   E. Bombieri, On exponential sums in finite fields, II, Invent. Math., 47(1978), 29-39.

[De]   P. Deligne, La conjecture de Weil, I, Publ. Math., IHES 43(1974), 273-307.

[Dw]   B. Dwork, On the rationality of the zeta function of an algebraic variety, Amer. J. Math., 82(1960), 631-648.

[ES]   J. -Y. Etesse and B. Le Stum, Fonctions L associées aux F-isocristaux surconvergents II, zéros et pôles unités, Invent. Math., 127(1997), no.1, 1-31.

[FHJ]  M. Fried, D. Haran and M. Jarden, Effective counting of the points of definable sets over finite fields, Israel J. Math., 85(1994), no. 1-3, 103-133.

[Gr]   A. Grothendieck, Formule de Lefschetz et rationalité des fonctions L, Séminaire Bourbaki, exposé 279, 1964/65.

[GA]   H. Gunji and D. Arron, On polynomial factorization over finite fields, Math. Comp., 36(1981), 281-287.

[GKS]  J. von zur Gathen, M. Karpinski and I. Shparlinski, Counting curves and their projections, Comput. Complexity, 6(1996/97), no.1, 64-99.

[Ka1]  N. Katz, Travaux de Dwork, Séminaire Bourbaki, Exp. 409(1971/72), Lecture Notes in Math. 317, 1973, 167-200.

[Ka2]  N. Katz, Une formule de congruence pour la fonction $\zeta$, in: Groupes de monodromie en géométrie algébrique (SGA 7 II), Exp. XXII, 401-438, Lecture Notes in Math. 340, Springer-Verlag, 1973

[Ko1]  N. Koblitz, Elliptic curve cryptosystems, Math. Comp., 48(1987), 203-209.

[Ko2]  N. Koblitz, Hyperelliptic cryptosystems, J. Cryptology, 1(1989), 139-150.

[LV]   T.C.Y. Lee and S.A. Vanstone, Subspaces and polynomial factorizations over finite fields, AAECC, 6(1995), 147-157.

[Mo]   P. Monsky, Formal cohomology III, Ann. Math., 93(1971), 315-343.

[Ni]   H. Niederreiter, Factoring polynomials over finite fields using differential equations and normal bases, Math. Comp., 62(1994), 819-830.

[Pi]   J. Pila, Frobenius maps of abelian varieties and finding roots of unity in finite fields, Math. Comp., 55(1990), 745-763.

[Po]   B. Poonen, Computational aspects of curves of genus at least 2, Lecture Notes in Computer Science, 1122 (1996), 283-306.

[Sc1]  R. Schoof, Elliptic curves over finite fields and the computation of square roots modulo $p$, Math. Comp., 44(1985), 483-494.

[Sc2]  R. Schoof, Counting points on elliptic curves over finite fields, Journal de Théorie des Nombres de Bordeaux, 7(1995), 219-254.

[TW1]  Y. Taguchi and D. Wan, $L$-functions of Drinfeld modules and $\varphi$-sheaves, J. Amer. Math. Soc., 9(1996), 755-781.

[TW2]  Y. Taguchi and D. Wan, Entireness of $L$-functions of $\varphi$-sheaves on affine complete intersections, J. Number Theory, 63(1997), no.1, 170-179.

[Ts]   M.A. Tsfasman, Algebraic geometry lattices and codes, Lecture Notes in Computer Science, 1122 (1996), 385-389.

[Wa1]  D. Wan, Hilbert sets and zeta functions over finite fields, J. Reine Angew. Math., 427(1992), 193-207.

[Wa2]  D. Wan, Meromorphic continuation of L-functions of $p$-adic representations, Ann. Math., 143(1996), 469-498.

DEPARTMENT OF MATHEMATICS, UNIVERSITY OF CALIFORNIA, IRVINE, CA 92697-3875
*E-mail address*: `dwan@math.uci.edu`